\documentclass[10pt,a4paper]{amsart}

\usepackage[english]{babel}

\usepackage{color}
\usepackage{amscd}
\usepackage{amssymb}
\usepackage{euscript}

\newcommand{\UnTrait}{
\bigskip \begin{center} \rule{6cm}{0.1mm} \end{center}}

\newtheorem{theorem}{Theorem}[subsection]

\newtheorem{proposition}[theorem]{Proposition}
\newtheorem{lemma}[theorem]{Lemma}
\newtheorem{corollary}[theorem]{Corollary}

\theoremstyle{remark}
\newtheorem{remark}[theorem]{Remark}
\newtheorem{example}[theorem]{Example}

\renewcommand{\mod}{\mathrm{mod}~}

\renewcommand{\Im}{\mathop{\mathrm{Im}}\nolimits}

\newcommand{\card}[1]{\lvert#1\rvert}
\newcommand{\Ch}[2]{\begin{bmatrix} #1 \\ #2 \end{bmatrix}}
\newcommand{\Id}{\textrm{Id}}
\newcommand{\idest}{\textit{i.~e.}}
\newcommand{\Inv}{\mathsf{Inv}}
\newcommand{\isom}{{\stackrel{\sim}{\longrightarrow}}}
\newcommand{\loccit}{\textit{loc.~cit.}}
\newcommand{\modtimes}{\mathrm{mod^{\times}}~}
\newcommand{\set}[2]{\left\{#1\,\mid\,#2\right\}}
\newcommand{\Sym}{\textrm{Sym}}

\DeclareMathOperator{\chr}{char}
\DeclareMathOperator{\Disc}{Disc}
\DeclareMathOperator{\Gal}{Gal}
\DeclareMathOperator{\GL}{GL}

\DeclareMathOperator{\Jac}{Jac}

\DeclareMathOperator{\Res}{Res}
\DeclareMathOperator{\SL}{SL}
\DeclareMathOperator{\Sp}{Sp}

\newcommand{\Am}{\mathsf{A}}
\newcommand{\Cm}{\mathsf{C}}
\newcommand{\Mm}{\mathsf{M}}

\newcommand{\Vm}{\mathsf{V}}
\newcommand{\Xm}{\mathsf{X}}
\newcommand{\Ym}{\mathsf{Y}}

\newcommand{\am}{\boldsymbol{\alpha}}
\newcommand{\epsm}{\boldsymbol{\varepsilon}}

\newcommand{\lm}{\boldsymbol{\lambda}}
\newcommand{\wm}{\boldsymbol{\omega}}

\newcommand{\CC}{\mathbb{C}}
\newcommand{\HH}{\mathbb{H}}
\newcommand{\NN}{\mathbb{N}}
\newcommand{\PP}{\mathbb{P}}
\newcommand{\QQ}{\mathbb{Q}}

\newcommand{\UU}{\mathbb{U}}
\newcommand{\ZZ}{\mathbb{Z}}

\newcommand{\bfI}{\mathbf{I}}
\newcommand{\bfM}{\mathbf{M}}
\newcommand{\bfR}{\mathbf{R}}
\newcommand{\bfS}{\mathbf{S}}
\newcommand{\bfT}{\mathbf{T}}

\begin{document}

\title[Jacobians among Abelian threefolds]{Jacobians among Abelian
threefolds: \\ a formula of Klein and a question of Serre}

\author{Gilles Lachaud}
\address{
Gilles Lachaud
\newline \indent
Institut de Math\'ematiques de Luminy
\newline \indent
Universit\'e Aix-Marseille - CNRS
\newline \indent
Luminy Case 907, 13288 Marseille Cedex 9 - FRANCE
\newline \indent}
\email{lachaud@iml.univ-mrs.fr}

\author{Christophe Ritzenthaler}
\address{
Christophe Ritzenthaler
\newline \indent
Institut de Math\'ematiques de Luminy
\newline \indent}
\email{ritzent@iml.univ-mrs.fr}

\author{Alexey Zykin}
\address{
Alexey Zykin
\newline \indent
Institut de Math\'ematiques de Luminy
\newline \indent
Mathematical Institute of the Russian Academy of Sciences
\newline \indent
Laboratoire Poncelet (UMI 2615)
\newline \indent}
\email{zykin@iml.univ-mrs.fr}

\date{\today}


\date{\today}

\maketitle


\begin{abstract}
Let $k$ be a field and $f$ be a Siegel modular form of weight $h \geq 0$
and genus $g>1$ over $k$. Using $f$, we define an invariant of the
$k$-isomorphism class of a principally polarized abelian variety $(A,a)/k$
of dimension $g$. Moreover when $(A,a)$ is the Jacobian of a smooth plane
curve, we show how to associate to $f$ a classical plane invariant. As
straightforward consequences of these constructions when $g=3$ and $k \subset \CC$ we obtain 
(i) a new proof of a formula of Klein linking the modular form $\chi_{18}$
to the square of the discriminant of plane quartics ; (ii) a proof that
one can decide when $(A,a)$ is a Jacobian over $k$ by looking whether the value
of $\chi_{18}$ at $(A,a)$ is a square in $k$. This answers a
question of J.-P.~Serre.  Finally, we study the possible generalizations of this approach for
$g>3$.
\end{abstract}

\section{Introduction}

\subsection{Torelli theorem}

Let $k$ be an algebraically closed field. If $X$ is a
(nonsingular irreducible projective) curve of genus $g$ over $k$, Torelli's
theorem states that the map $X \mapsto (\Jac X,j),$ associating to $X$ its
Jacobian together with the canonical polarization $j$, is injective. The
determination of the image of this map is a long time studied question.

When $g = 3$, the moduli space $\Am_g$ of principally polarized abelian
varieties of dimension $g$ and the moduli space $\Mm_g$ of nonsingular
algebraic curves of genus $g$ are both of dimension $3g-3=g (g+1)/2=6$.
According to Hoyt \cite{Hoyt} and Oort and Ueno \cite{Ueno}, the image of
$\Mm_3$ is exactly the space of indecomposable principally polarized
abelian threefolds. Moreover if $k=\CC$, Igusa \cite{igusa2}
characterized the locus of decomposable abelian threefolds and that of
hyperelliptic Jacobians making use of two particular modular forms
$\Sigma_{140}$ and $\chi_{18}$ on the Siegel upper half space of degree
$3$.

Assume now that $k$ is any field and $g \geq 1$. J.-P. Serre noticed in
\cite{lauter} that a precise form of Torelli's theorem reveals a mysterious
obstruction for a geometric Jacobian to be a Jacobian over $k$. More
precisely, he proved the following:

\begin{theorem}
\label{twist}
Let $(A,a)$ be a principally polarized abelian variety  of dimension $g>0$
over $k$, and assume that $(A,a)$ is isomorphic over $\overline{k}$ to
the Jacobian of a curve $X_{0}$ of genus $g$ defined over $\overline{k}$.
The following alternative holds :
\begin{enumerate}
\item
\label{TWI1}
If $X_{0}$ is hyperelliptic, there is a curve $X/k$ isomorphic to
$X_{0}$ over $\overline{k}$ such that $(A,a)$ is $k$-isomorphic
to $(\Jac X,j)$.
\item
\label{TWI2}
If $X_{0}$ is not hyperelliptic, there is  a curve $X/k$ isomorphic to
$X_{0}$ over $\overline{k}$, and a quadratic character
$$
\begin{CD}
\varepsilon : \Gal(k_{\textrm{sep}}/k) & @>>> &
\{\pm 1\}
\end{CD}
$$
such that the twisted abelian variety
$(A,a)_{\varepsilon}$ is $k$-isomorphic to $(\Jac X,j)$. The character
$\varepsilon$ is trivial if and only if $(A, a)$ is $k$-isomorphic to a
Jacobian.
\end{enumerate}
\end{theorem}

Thus, only case \eqref{TWI1} occurs if $g = 1$ or $g = 2$, with all curves
being elliptic or hyperelliptic.

\subsection{Curves of genus $3$}

Assume now $k \subset \CC$ and $g = 3$.
Let there be given an indecomposable  principally polarized abelian
threefold $(A, a)$ defined over $k$. In a letter to J. Top
\cite{SerreTop}, J.-P. Serre asked a twofold question:

\begin{itshape}
\begin{itemize}
\item How to decide, knowing only $(A, a)$, that $X$ is hyperelliptic ?
\item If $X$ is not hyperelliptic, how to find the quadratic character $\varepsilon$ ?
\end{itemize}
\end{itshape}

Moreover, he suggested a strategy in order to compute the twisting factor $\varepsilon$. This strategy is based on a formula of Klein \cite{klein} relating the modular form $\chi_{18}$ (in the notation of Igusa), to the square of the discriminant of plane quartics, see Th.\ref{kleinformula} for a more precise formulation. In \cite{SAGA}, two of the authors justified Serre's strategy for a three dimensional family of abelian varieties and in particular determined the absolute constant involved in Klein's formula.

In this article we prove that Serre's strategy can be applied to any abelian threefolds. More
precisely, we take a broader point of view.

\begin{enumerate}
\item \label{point1}
We look at the action of $\overline{k}$-isomorphisms
on Siegel modular forms defined over $k$ and we define
 invariants of  $k$-isomorphism classes of abelian varieties over $k$.
\item \label{point2}
We link Siegel modular forms, Teichm\"{u}ller modular forms and
invariants. Then we derive a proof of Klein's formula based on moduli spaces.
\end{enumerate}

Once these two goals achieved, Serre's strategy can be rephrased as finding a Siegel modular form whose locus has a good
multiplicity on the Jacobian locus and then using point (\ref{point1}) to distinguish between
Jacobians and their twists. For $g = 3$, the form $\chi_{18}$
fulfills the criterion as can be seen thanks to Klein's formula.
On the other hand, we show that this is no longer
the case for $\chi_{h}$ when $g>3$.
We would like to point out that we do not actually
need Klein's formula to prove Serre's strategy.
Indeed we do not need to go the full way from
Siegel modular form to invariants and could instead use a formula due to
Ichikawa relating $\chi_{18}$ to the square of a Teichm\"{u}ller modular form (see Rem. \ref{mu9serre}).
However we think that the connection between Siegel modular forms and invariants
is interesting enough in its own, besides the fact that it gives a new proof of
Klein's formula.\\

The paper is organized as follows. In \S \ref{part1}, we review the
necessary elements from the theory of Siegel and Teichm\"{u}ller modular
forms. Only \S \ref{partaction} is original: we introduce the action of
isomorphisms and see how the action of twists is reflected on the values of
modular forms. In \S \ref{part2}, we link modular forms and certain
invariants of ternary forms. Finally in \S \ref{part3} we deal with the
case $g=3$. We give first a proof of Klein's formula and then we justify the
validity of Serre's strategy. Finally we explain the reasons behind the
failure of the obvious generalization of the theory in higher dimensions and state some natural questions.\\

{\bf Acknowledgements.} We would like to thank J.-P. Serre and S. Meagher for fruitful discussions and Y. F. Bilu and X. Xarles for their help in the final part of Sec.\ref{beyond}.

\section{Siegel and Teichm\"{u}ller  modular forms}
\label{part1}

\subsection{Geometric Siegel modular forms}

The references are \cite{Chai}, \cite{DM}, \cite{FC}, \cite{GVDG}. Let $g
> 1$ and $n > 0$ be two integers and $\Am_{g,n}$ be the moduli stack of
principally polarized abelian schemes of relative dimension $g$ with
symplectic level $n$ structure. Let $\pi : \Vm_{g,n} \longrightarrow
\Am_{g,n}$ be the universal abelian scheme, fitted with the zero section
$\varepsilon:
\Am_{g,n}
\longrightarrow \Vm_{g,n}$, and
$$
\begin{CD}
\pi_{*} \Omega^{1}_{\Vm_{g,n}/\Am_{g,n}} =
\varepsilon^{*}  \Omega^{1}_{\Vm_{g,n} / \Am_{g,n}} & @>>> & \Am_{g,n}
\end{CD}
$$
the rank $g$ bundle induced by the  relative regular differential forms of
degree one on $\Vm_{g,n}$ over $\Am_{g,n}$. The relative canonical bundle
over $\Am_{g,n}$ is the line bundle
$$
\wm =\wedge^{g} \varepsilon^{*} \Omega^{1}_{\Vm_{g,n} / \Am_{g,n}}.
$$
For a projective nonsingular variety $X$ defined over a
field $k$, we denote by
$$\Omega^{1}_{k}[X] = H^{0}(X,\Omega^{1}_{X} \otimes k)$$
the finite dimensional $k$-vector space of regular differential forms on
$X$ defined over $k$. Hence, the fibre of the
bundle $\Omega^{1}_{\Vm_{g,n} / \Am_{g,n}}$ over $A \in \Am_{g,n}(k)$ is equal to
$\Omega^{1}_{k}[A]$, and the fibre of $\wm$ is the
one-dimensional vector space
$$\wm[A]= \wedge^{g} \Omega^{1}_{k}[A].$$
We put $\Am_g=\Am_{g,1}$ and $\Vm_{g} = \Vm_{g,1}$. Let $R$ be a
commutative ring and $h$ a positive integer. A \emph{geometric
Siegel modular form} of genus $g$ and weight $h$ over  $R$ is an element
of the
$R$-module
$$
\bfS_{g,h}(R) = \Gamma(\Am_{g} \otimes R,\wm^{\otimes h}).
$$
Note that for any $n \geq 1$, we have an isomorphism
$$\Am_{g} \simeq \Am_{g,n}/\Sp_{2g}(\ZZ/n\ZZ).$$
If $n \geq 3$, as shown in \cite{GIT}, from the rigidity lemma of
Serre \cite{SerreJacobi} we can deduce that the moduli space $\Am_{g,n}$ can be
represented by a smooth scheme over $\ZZ[\zeta_{n},1/n]$. Hence, for
any algebra $R$ over $\ZZ[\zeta_{n},1/n]$, the module $\bfS_{g,h}(R)$ is
the submodule of
$$\Gamma(\Am_{g,n} \otimes_{\ZZ[\zeta_n,1/n]} R,\wm^{\otimes h})$$
consisting of the elements invariant under $\Sp_{2g}(\ZZ/n\ZZ)$.

Assume now that $R=k$ is a field. If $f \in \mathbf{S}_{g,h}(k)$, $A$ is a
p.p.a.v.  of dimension $g$ defined over $k$ and $\alpha$ is a basis of
$\boldsymbol{\omega}_k[A]$, define

\begin{equation}
\label{DSMF}
f(A,\alpha) = f(A)/ \alpha^{\otimes h}.
\end{equation}

In this way such a modular form defines a rule which assigns the element
$f(A,\alpha)\in k$ to every such pair $(A,\alpha)$ and such that:

\begin{enumerate}
\item
$f(A, \lambda \alpha) = \lambda^{- h} f(A,\alpha)$ for any $\lambda \in
k^{\times}$.
\item
$f(A,\alpha)$ depends only on the $\overline{k}$-isomorphism class of the pair
$(A,\alpha)$.
\end{enumerate}

Conversely, such a rule defines a unique $f \in \mathbf{S}_{g,h}(k)$.
This definition is a straightforward generalization of that of Deligne-Serre \cite{DS} and
Katz \cite {Katz} if $g = 1$.

\subsection{Complex uniformisation}

Assume $R = \CC$. Let
$$
\HH_{g} =
\set{\tau \in \mathbf{M}_{g}(\CC)}{^{t}{\tau} = \tau, \ \Im \tau > 0}
$$
be the Siegel upper half space of genus $g$ and $\Gamma = \Sp_{2g}(\ZZ)$.
As explained in \cite[\S 2]{Chai}, The complex orbifold $\Am_g(\CC)$ can be
expressed as the quotient
$\Gamma\backslash\HH_g$ where
$\Gamma$ acts by
$$
M.\tau = (a \tau + b)\cdot(c \tau + d)^{-1}
\quad \text{if} \quad
M = \begin{pmatrix} a & b \\ c & d \end{pmatrix} \in \Gamma.
$$
The group $\ZZ^{2 g}$ acts on $\HH_{g} \times \CC^{g}$ by
$$
v.(\tau, z) = (\tau, z + \tau m + n)
\quad \text{if} \quad
v = \begin{pmatrix} m \\ n \end{pmatrix} \in \ZZ^{2 g}.
$$
If  $\UU_{g} = \ZZ^{2 g} \backslash (\HH_{g} \times \CC^{g})$, the
projection
$$
\begin{CD}
\pi : \UU_{g} & @>>> & \HH_{g}
\end{CD}
$$
defines a universal principally polarized abelian variety with fibres
$$A_{\tau} = \pi^{-1}(\tau)=\CC^{g}/(\ZZ^{g} + \tau \ZZ^{g}).$$
Let $j(M,\tau) = c  \tau + d$ and define the action of $\Gamma$ on
$\HH_{g}
\times
\CC^{g}$ by
$$
M.(\tau,(z_{1}, \dots, z_{g})) =
(M.\tau, {^{t}j(M,\tau)}^{-1} \cdot (z_{1}, \dots, z_{g}))
\quad \text{if} \ M \in \Gamma.
$$
The map ${^{t}j(M,\tau)}^{-1} : \CC^{g} \rightarrow \CC^{g}$ induces an
isomorphism:
$$
\begin{CD}
\varphi_{M} : A_{\tau} & @>>> & A_{M.\tau}. \\
\end{CD}
$$
Hence, $\Vm_{g}(\CC) \simeq \Gamma \backslash \UU_{g}$ and the following
diagram is commutative:
$$
\begin{CD}
\Gamma \backslash \UU_{g} & @>{\sim}>> & \Vm_{g}(\CC) \\
@V{\pi}VV & & @V{\pi}VV \\
\Gamma \backslash \HH_{g} & @>{\sim}>> & \Am_{g}(\CC)
\end{CD}
$$
As in \cite[p. 141]{FC}, let
$$
\zeta =
\frac{dq_{1}}{q_{1}}\wedge \dots \wedge \frac{dq_{g}}{q_{g}} =
(2 i \pi)^{g}  dz_{1} \wedge \cdots \wedge dz_{g} \in
\Gamma(\HH_{g},\wm)
$$
with  $(z_{i}, \dots, z_{g}) \in \CC^{g}$ and $(q_{i}, \dots, q_{g}) =
(e^{2 i \pi z_{1}}, \dots e^{2 i \pi z_{g}})$. This section of the canonical
bundle is a basis of
$\wm[A_{\tau}]$ for all  $\tau \in \HH_{g}$
and the relative canonical bundle of $\UU_{g}/\HH_{g}$ is trivialized by
$\zeta$ :
$$
\wm_{\UU_{g}/\HH_{g}} =
\wedge^{g} \Omega^{1}_{\UU_{g}/\HH_{g}} \simeq
\HH_{g} \times \CC \cdot \zeta.
$$
The group $\Gamma$ acts on $\wm_{\UU_{g}/\HH_{g}}$ by
$$
M.(\tau, \zeta) =
(M.\tau, \det j(M,\tau) \cdot \zeta)
\quad \text{if} \ M \in \Gamma,
$$
in such a way that
$$
\varphi_{M}^{*}(\zeta_{M.\tau}) = \det j(M,\tau)^{- 1}\zeta_{\tau}.
$$
Thus, a geometric Siegel modular form $f$ of weight $h$ becomes an
expression
$$f(A_{\tau}) = \widetilde{f}(\tau) \cdot \zeta^{\otimes h},$$
where $\widetilde{f}$ belongs to the well-known vector space $\mathbf{R}_{g,h}(\CC)$
of \emph{analytic Siegel modular forms} of weight $h$ on $\HH_{g}$, consisting of complex holomorphic functions $\phi(\tau)$ on
$\HH_{g}$ satisfying
$$
\phi(M.\tau) = \det j(M.\tau)^{h} \phi(\tau)
$$
for any $M \in \Sp_{2g}(\ZZ)$. Note that by Koecher principle \cite[p.
11]{GVDG}, the condition of holomorphy at $\infty$ is automatically
satisfied since $g > 1$. The converse is also true \cite[p. 141]{FC}:

\begin{proposition}
\label{HoloGeom}
If $f \in \bfS_{g,h}(\CC)$ and $\tau \in \HH_{g}$, let
$$
\widetilde{f}(\tau) = f(A_{\tau})/\zeta^{\otimes h} = (2i\pi)^{- g h}
f(A_{\tau})/
(dz_{1} \wedge \cdots \wedge dz_{g})^{\otimes h}.
$$
Then the map $f \mapsto \widetilde{f}$ is an isomorphism
$\bfS_{g,h}(\CC) \isom \bfR_{g,h}(\CC).$
\qed
\end{proposition}

\subsection{Teichm\"{u}ller modular forms}
\label{sec_Teichmuller}

Let $g > 1$ and $n > 0$ be positive integers and let $\Mm_{g,n}$ denote the
moduli stack of smooth and proper curves of genus $g$ with symplectic level
$n$ structure \cite{DM}. Let $\pi : \Cm_{g,n} \longrightarrow \Mm_{g,n}$ be
the universal curve, and let $\lm$ be the invertible sheaf associated to the
\emph{Hodge bundle}, namely
$$
\lm = \wedge^{g} \pi_{*}\Omega^{1}_{\Cm_{g,n} / \Mm_{g,n}}.
$$
For an algebraically closed field $k$ the fibre over $C \in \Mm_{g,n}(k)$
is the one dimensional vector space $\lm[C] = \wedge^{g} \Omega^{1}_{k}[C]$.

Let $R$ be a commutative ring and $h$ a positive integer. A
\emph{Teichm\"{u}ller modular form}\ of genus $g$ and weight $h$ over  $R$
is an element of
$$\bfT_{g,h}(R) = \Gamma(\Mm_{g} \otimes R, \lm^{\otimes h}).
$$
These forms have been thoroughly studied by Ichikawa \cite{Ichi1},
\cite{Ichi2}, \cite{Ichi3}, \cite{Ichi4}. As in the case of the moduli
space of abelian varieties, for any $n \geq 1$ we have
$$\Mm_g \simeq \Mm_{g,n}/\Sp_{2g}(\ZZ/n\ZZ),$$
and $\Mm_{g,n}$ can be represented by a smooth scheme over
$\ZZ[\zeta_{n},1/n]$ if $n \geq 3$. Then, for any algebra $R$ over
$\ZZ[\zeta_n,1/n]$, the module $\bfT_{g,h}(R)$ is the submodule of
$$
\Gamma(\Mm_{g,n} \otimes_{\ZZ[\zeta_n,1/n]} R,\lm^{\otimes h})
$$
invariant under $\Sp_{2g}(\ZZ/n\ZZ)$.

Let $C/k$ be a genus $g$ curve. Let $\lambda_{1}, \dots, \lambda_{g}$ be a
basis of $\Omega^{1}_{k}[C]$ and $\lambda = \lambda_{1} \wedge \dots \wedge
\lambda_{g}$ a basis of $\lm [C]$. As for Siegel modular forms in (\ref{DSMF}), for a Teichm\"{u}ller modular form $f\in \bfT_{g,h}(k)$ we
define
$$f(C,\lambda)=f(C)/\lambda^{\otimes h} \in k.$$

Ichikawa proves the following proposition:

\begin{proposition}
\label{Tau}
The Torelli map $\theta : \Mm_{g} \longrightarrow\Am_{g}$, associating to a
curve $C$ its Jacobian $\Jac C$ with the canonical polarization $j$,
satisfies  $\theta^{*}\wm = \lm$, and induces for any commutative ring $R$
a linear map
$$
\begin{CD}
\theta^{*} :
\bfS_{g,h}(R) =
\Gamma(\Am_{g} \otimes R, \wm^{\otimes h}) & @>>> &
\bfT_{g,h}(R) =
\Gamma(\Mm_{g} \otimes R, \lm^{\otimes h}),
\end{CD}
$$
such that $[\theta^{*}f](C) = \theta^{*}[f(\Jac C)]$.
Fixing a basis $\lambda$ of $\lm [C]$, this is
$$f(\Jac C,\alpha) = [\theta^{*}f](C,\lambda) \qquad \text{if} \ \theta^{*}\alpha = \lambda.$$
\qed
\end{proposition}

\subsection{Action of isomorphisms}
\label{partaction}
Suppose $\phi : (A',a') \longrightarrow (A,a)$ is a
$\overline{k}$-isomorphism of principally polarized abelian varieties, then
by definition
$$f(A,\alpha) = f(A',\beta)$$
where $\beta_{i} = \phi^{*}(\alpha_{i})$ is a basis of
${\Omega^{1}}_{\overline{k}}[A']$ and $\beta = \beta_{1} \wedge \dots \wedge
\beta_{g} \in \wm[A']$. If $\alpha'_{1}, \dots, \alpha'_{g}$ is another
basis of ${\Omega^{1}}_{\overline{k}}[A']$ and $\alpha' = \alpha'_{1} \wedge
\dots \wedge \alpha'_{g}$, we denote by $M_{\phi}\in \GL_{g}(\overline{k})$
the matrix of the basis $(\beta_{i})$ in the basis $(\alpha'_{i})$. We can easily see that

\begin{proposition}
\label{det}
In the above notation,
$$f(A,\alpha) = \det(M_{\phi})^h \cdot f(A',\alpha'). \rlap \qed$$
\end{proposition}

First of all, from this formula applied to the action of $-1,$ we deduce
that, if $k$ is a field of characteristic different from $2$, then
$\bfS_{g,h}(k) = \{0\}$ if $gh$ is odd. From now on we assume that
$gh$ is even and $\chr k \neq 2$.

\begin{corollary}
\label{ftwist}
Let $(A,a)$ be a principally polarized abelian variety of dimension $g$
defined over $k$ and $f \in \bfS_{g,h}(k)$. Let $\alpha_{1}, \dots,
\alpha_{g}$ be a basis of $\Omega^{1}_{k}[A]$, and put $\alpha
= \alpha_{1} \wedge \dots \wedge \alpha_{g} \in \wm[A]$. Then the quantity
$$\bar{f}(A) = f(A,\alpha) \ \modtimes k^{\times h} \in k/k^{\times
h}$$ does not depend on the choice of the basis of $\Omega^{1}_{k}[A]$. In
particular $\bar{f}(A)$ is an invariant of the $k$-isomorphism class of $A$.
\qed
\end{corollary}

\begin{corollary}
\label{cortwist}
Assume that $g$ is odd. Let $f \in \bfS_{g,h}(k)$ and $\phi : A'
\longrightarrow A$ a non trivial  quadratic twist. If $\bar{f}(A) \neq
0$ then $\bar{f}(A)$ and
$\bar{f}(A')$ do not belong to the same class in $k^{\times} /
k^{\times 2}$.
\end{corollary}

\begin{proof}
Assume that $\phi$ is given by the quadratic character $\varepsilon$ of
$\Gal(\overline{k}/k)$. Then
$$
d^{\sigma} = \varepsilon(\sigma)^g \cdot  d, \text{ where} \quad
d = \det(M_{\phi}) \in \overline{k},\quad \sigma \in \Gal(\overline{k}/k).
$$
Assume that $g$ is odd. Then by our assumption $h$ is even, and
$d^2 = \varepsilon(\sigma)^g d d^{\sigma} \in k$. But $d \notin k$ since
there exists $\sigma$ such that $\varepsilon(\sigma) = -1$. Using
Prop.\ref{det} we find that
$$f(A,\alpha)= (d^2)^{h/2} f(A',\alpha').$$
Since $d^2$ is not a square in $k$, if $\bar{f}(A) \neq 0$ then
$\bar{f}(A)$ and $\bar{f}(A')$ belong to two different classes
in $k^{h/2}/k^{\times h} \simeq k/k^{\times 2}$.
\end{proof}

Let now $(A,a)$ be a principally polarized abelian variety of dimension $g$
defined over $\CC$.
Let $\omega_{1}, \dots, \omega_{g}$ be a basis of $\Omega^{1}_{\CC}[A]$
and $\omega = \omega_{1} \wedge \dots \wedge \omega_{g} \in \wm[A]$. Let
$\gamma_{1}, \dots \gamma_{2g}$ be a symplectic basis (for the polarization
$a$). The period matrix
$$
\Omega = [\Omega_{1} \ \Omega_{2}] =
\begin{pmatrix}
\int_{\gamma_{1}}\omega_{1} & \cdots & \int_{\gamma_{2g}}\omega_{1}\\
\vdots &  & \vdots \\
\int_{\gamma_{1}}\omega_{g} & \cdots & \int_{\gamma_{2g}}\omega_{g}\\
\end{pmatrix}
$$
belongs to the set $\EuScript{R}_{g} \subset \bfM_{g,2g}(\CC)$ of
Riemann matrices, and $\tau = \Omega_{2}^{- 1} \Omega_{1} \in \HH_{g}$.

\begin{proposition}
\label{Riemann}
In the above notation,
$$
f(A,\omega)=(2i \pi)^{gh}  \frac{\widetilde{f}(\tau)}{\det \Omega_{2}^{h}}.
$$
\end{proposition}

\begin{proof}
The
abelian variety $A$ is isomorphic to $A_{\Omega} = \CC^{g} / \Omega
\ZZ^{2g}$ and $\omega \in \wm[A]$ maps to $\xi = dz_{1} \wedge \cdots \wedge
dz_{g} \in \wm[A_{\Omega}]$ under this isomorphism. The linear map  $z \mapsto \Omega_{2}^{-1}z=z'$
induces the isomorphism
$$
\begin{CD}
\varphi :  A_{\Omega} & @>>> & A_{\tau} = \CC^{g}/(\ZZ^{g}+\tau\ZZ^{g}).
\end{CD}
$$
Let us denote $\xi'=dz_1' \wedge \cdots \wedge dz'_g=(2i\pi)^{-g} \zeta$ in $\wm[A_{\tau}]$. Thus, using Prop.\ref{det}, Equation \eqref{DSMF} and Prop.\ref{HoloGeom},
we obtain
\begin{multline*}
f(A,\omega) = f(A_\Omega,\xi) =
\det \Omega_{2}^{- h} f(A_{\tau},\xi' )
\\
= \det \Omega_{2}^{- h} f(A_{\tau})/ {\xi'}^{\otimes h} =
(2i \pi)^{gh} \det \Omega_{2}^{- h} f(\tau)/ \zeta^{\otimes h} =
(2i \pi)^{gh} \frac{\widetilde{f}(\tau)}{\det \Omega_{2}^{h}},
\end{multline*}
from which the proposition follows.
\end{proof}

\section{Invariants and modular forms}
\label{part2}
In this section $k$ is an algebraically closed field of characteristic
different from $2$.

\subsection{Invariants}
\label{Invariants}

We review some classical invariant theory. Let $E$ be a vector space of
dimension $n$ over $k$. The left regular representation $\rho$ of $\GL(E)$
on the vector space $\Xm_{d} = \Sym^{d}(E^{*})$ of homogeneous
polynomials of degree $d$ on $E$ is given by
$$
\rho(u): F(x) \mapsto
(u \cdot F)(x) = F(u^{-1}x)
$$
for $u \in \GL(E)$, $F \in \Xm_{d}$ and $x \in E$. If $U$ is an open subset
of $\Xm_{d}$ stable under $\rho$, we still denote by $\rho$ the left
regular representation of $\GL(E)$ on the $k$-algebra $\mathcal{O}(U)$ of
regular functions on $U$, in such a way that
$$
\rho(u): \Phi(F) \mapsto (u \cdot \Phi)(F) = \Phi(u^{-1} \cdot F),
$$
if $u \in \GL(E)$, $\Phi \in \mathcal{O}(U)$ and $F \in U$.  If $h \in \ZZ$,
we denote by $\mathcal{O}_{h}(U)$ the subspace of homogeneous elements
of degree $h$, satisfying $\Phi(\lambda F) = \lambda^{h} \Phi(F)$ for
$\lambda \in k^{\times}$ and $F \in U$. The subspaces $\mathcal{O}_{h}(U)$
are stable under $\rho$. An element $\Phi \in
\mathcal{O}_{h}(U)$ is an \emph{invariant of degree $h$ on $U$} if
$$
u \cdot \Phi = \Phi  \quad \text{for every} \ u \in \SL(E),
$$
and we denote by $\Inv_{h}(U)$ the subspace of $\mathcal{O}_{h}(U)$ of
invariants of degree $h$ on $U$. If $\Inv_{h}(U) \neq \{0\}$, then $hd
\equiv 0 (\mod n)$, since the group $\boldsymbol{\mu}_{n}$ of $n$-th roots
of unity is in the kernel of $\rho$. Hence, if $\Phi \in \mathcal{O}(U)$,
and if $w$ and $n$ are two integers such that $h d = n w$, the following
statements are equivalent:
\begin{enumerate}
\item
$\Phi \in \Inv_{h}(U)$;
\item
$u \cdot \Phi = (\det u)^{-w} \Phi \quad \text{for every} \ u \in \GL(E).$
\end{enumerate}
If these conditions are satisfied, we call $w$ the \emph{weight} of $\Phi$.

The \emph{multivariate resultant} $\Res(f_{1},\dots,f_{n})$ of $n$
forms $f_{1}, \dots f_{n}$ in $n$ variables with coefficients in $k$ is an
irreducible polynomial in the coefficients of $f_{1}, \dots f_{n}$ which
vanishes whenever $f_{1}, \dots f_{n}$ have a common non-zero root. One
requires that the resultant is irreducible over $\ZZ$, \idest{} it has
coefficients in $\ZZ$ with greatest common divisor equal to $1$, and
moreover
$$\Res(x_{1}^{d_{1}},\dots,x_{n}^{d_{n}}) = 1$$
for any $(d_1,\ldots,d_n) \in \NN^n$. The resultant exists and is
unique. Now, let $F \in \Xm_{d}$, and denote $q_{1}, \dots, q_{n}$ the partial
derivatives of $F$. The \emph{discriminant} of $F$ is
$$
\Disc F = c_{n,d}^{- 1} \, \Res(q_{1},\dots,q_{n}), \quad \text{with} \quad
c_{n,d} = d^{((d - 1)^{n} - (- 1)^{n})/d},
$$
the coefficient $c_{n,d}$ being chosen according to \cite{SerreTop}. Hence,
the projective hypersurface which is the zero locus of $F \in \Xm_{d}$ is
nonsingular if and only if $\Disc F \neq 0$. The discriminant is an
irreducible polynomial in the coefficients of $F$, see for instance
\cite[Chap. 9, Ex. 1.6(a)]{GPZ}. From now on we restrict ourselves to the
case $n = 3$, \idest{} we consider invariants of ternary forms of degree
$d$, and summarize the results that we shall need.

\begin{proposition}
\label{discriminant}
If $F \in \Xm_{d}$ is a ternary form, the discriminant
$$\Disc F = d^{- (d - 1)(d - 2) - 1} \cdot \Res(q_{1},q_{2},q_{3})$$
where $q_{1}, q_{2}, q_{3}$ are the partial derivatives of $F$, is given by
an irreducible polynomial over $\ZZ$ in the coefficients of
$F$, and vanishes if and only if the plane curve $C_{F}$ defined by $F$
is singular. The discriminant is an invariant of
$\Xm_{d}$ of degree $3(d - 1)^{2}$ and weight $d(d - 1)^{2}$.
\qed
\end{proposition}

We refer to \cite[p. 118]{GPZ} and \cite{SAGA} for a beautiful explicit
formula for the discriminant, found by Sylvester.

\begin{example}[Ciani quartics]
\label{CianiQuartics}
We recall some results whose proofs are given in \cite{SAGA}. Let
$\Sym_{3}(k)$ be the vector space of symmetric matrices of size
$3$ with coefficients in $k$, and
$$
G_{m}(x, y, z) = {^{t}\!v}.m.v, \quad v = (x ,y ,z),
$$
the quadratic form associated to $m \in \Sym_{3}(k)$. Then
$$
F_{m}(x, y, z) = G_{m}(x^{2}, y^{2}, z^{2})
$$
is a ternary quartic, and the map $m \mapsto F_{m}$ is an
isomorphism of $\Sym_{3}(k)$ to the subspace of $F \in \Xm_{4}$ which are
invariant under the three involutions
$$
\sigma_{1}(x,y,z) = (-x,y,z), \quad \sigma_{2}(x,y,z) = (x,-y,z), \quad
\sigma_{3}(x,y,z) = (x,y,-z).
$$
If
$$
m =
\begin{bmatrix}
a_{1} & b_{3} & b_{2} \\
b_{3} & a_{2} & b_{1} \\
b_{2} & b_{1} & a_{3}
\end{bmatrix}
\in \Sym_{3}(k),
$$
then
$$
F_{m}(x, y, z) = a_{1} x^{4} + a_{2} y^{4} + a_{3} z^{4} + 2 (b_{1}
y^{2} z^{2} + b_{2} x^{2} z^{2} + b_{3} x^{2} y^{2}).
$$
For $1 \leq i \leq 3$, let $c_{i} = a_{j} a_{k} - b_{i}^2$ the cofactor
of $a_{i}$. Then
$$\Disc F_{m} = 2^{40} \, a_{1}\,a_{2}\,a_{3} \, (c_{1}\,c_{2}\,c_{3})^{2} \,
\det(m)^4.
$$
Note that the discrepancy between the powers of $2$ here and in \cite[Prop.2.2.1]{SAGA}
comes from the  normalization by $c_{n,d}$.
\end{example}

\subsection{Geometric invariants for nonsingular plane curves}
\label{Quartics}

Let $E$ be a vector space of dimension $3$ over $k$ and $G = \GL(E)$. The
\emph{universal curve} over the affine space $\Xm_{d}=\Sym^d(E)$ is the
variety
$$
\Ym_{d} =\set{(F,x) \in \Xm_{d} \times \PP^{2}}{F(x) = 0}.
$$
The \emph{nonsingular locus} of $X_{d}$ is the principal open set
$$
\Xm_{d}^{0} = (\Xm_{d})_{\Disc} = \set{F \in \Xm_{d}}{\Disc(F) \neq 0}.
$$

If $\Ym_{d}^{0}$ is the universal curve restricted to the
nonsingular locus, the projection is a smooth surjective $k$-morphism
$$
\begin{CD}
\Ym_{d}^{0} & @>>> & \Xm_{d}^{0}
\end{CD}
$$
whose fibre over $F$ is the non singular plane curve $C_F$.

We recall the classical way to write down an explicit $k$-basis of
$\Omega^{1}[C_{F}] =H^0(C_F,\Omega^1)$ for $F \in \Xm_{d}^{0}(k)$ (see
\cite[p. 630]{bk}). Let
$$
\eta_{1} = \frac{f(x_{2}dx_{3} - x_{3} dx_{2})}{q_{1}}, \quad
\eta_{2} = \frac{f(x_{3}dx_{1} - x_{1} dx_{3})}{q_{2}}, \quad
\eta_{3} = \frac{f(x_{1}dx_{2} - x_{2} dx_{1})}{q_{3}},
$$
where $q_{1}, q_{2}, q_{3}$ are the partial derivatives of $F$, and where
$f$ belongs to the space $\Xm_{d - 3}$ of ternary forms of degree $d
- 3$. The forms $\eta_{i}$ glue together and define a regular
differential form $\eta_{f}(F) \in \Omega^{1}[C_{F}]$. Since $\dim \,
\Xm_{d - 3} = (d - 1)(d - 2)/2 = g$, the linear map
$f \mapsto \eta_{f}(F)$ defines an isomorphism
$$
\begin{CD}
\Xm_{d - 3} & @>{\sim}>> & \Omega^{1}[C_{F}].
\end{CD}
$$
\newcommand{\relX}{\Ym_{d}^{0}/\Xm_{d}^{0}}\noindent
This implies that the sections
$\eta_{f} \in \Gamma(\Xm_{d}^{0}, \Omega^{1}_{\relX})$
lead to a trivialization
$$
\begin{CD}
\Xm_{d}^{0} \times \Xm_{d - 3} & @>{\sim}>> & \Omega^{1}_{\relX}.
\end{CD}
$$
An element $u \in G$ acts on $\Ym_{d}$ by
$$u \cdot (F,x) = (u \cdot F,u x),$$
and the projection
$\Ym_{d}^{0} \longrightarrow \Xm_{d}^{0}$
is $G$-equivariant.

We denote $\eta_{1}, \dots, \eta_{g}$ the sequence of sections obtained by
substituting for $f$ in $\eta_{f}$ the $g$ members of the canonical basis of
$\Xm_{d - 3}$, enumerated according to the lexicographic order, the
\emph{classical basis} of $\Gamma(\Xm_{d}^{0}, \Omega^{1}_{\relX})$. The
section
\begin{equation*}
\eta = \eta_{1} \wedge \dots \wedge \eta_{g}
\end{equation*}
is a basis of the one-dimensional space $\Gamma(\Xm_{d}^{0}, \am)$, where
$$
\am =\wedge^{g} \pi_{*} \Omega^{1}_{\relX},
$$
is the Hodge bundle of the universal curve over $\Xm_{d}^{0}$.
For every $F \in \Xm_{d}^{0}$, an element $u \in G$
induces by restriction an isomorphism
$$
\begin{CD}
\varphi_{u} : C_{ F } & @>>> & C_{u \cdot F},
\end{CD}
$$
which itself defines a linear automorphism $\varphi^{*}_{u}$ of $\am$.\\

For any $h \in \ZZ$, we denote by $\Gamma(\Xm_{d}^{0}, \am^{\otimes
h})^{G}$ the subspace of sections $s \in \Gamma(\Xm_{d}^{0},
\am^{\otimes h})$ such that
$$
\varphi_{u}^{*}(s) = s \quad \text{for every} \ u \in G.
$$
If $\alpha \in  \Gamma(\Xm_{d}^{0}, \am)$ and $F \in \Xm_d^0$, we define,
in the same way as in Equation \eqref{DSMF},
$$s(F,\alpha)=s(F)/\alpha^{\otimes h}.$$
Hence, $s \in \Gamma(\Xm_{d}^{0}, \am^{\otimes h})^{G}$
if and only if for all $u \in G$ and $F \in X_d^0$, one has
$$(\varphi_u^*s)(F,\alpha)=s(F,\alpha).$$

\begin{proposition}
\label{Serre}
The section $\eta \in \Gamma(\Xm_{d}^{0},
\am)$ satisfies the following properties.
\begin{enumerate}
\item
\label{SER2}
If $u \in G$, then
$$
\varphi_u^* \eta=\det(u)^{ w_0}\eta,\quad \text{with}
\ w_0 = \binom{d}{3} = \frac{d g}{3} \, \in \NN.
$$
\item
\label{SER3}
Let $h\geq 0$ be an integer. The linear map
$$
\Phi \mapsto \tau(\Phi) = \Phi \cdot \eta^{\otimes h}
$$
is an isomorphism
$$
\begin{CD}
\tau :
\Inv_{gh}(\Xm_{d}^{0})& @>{\sim}>> & \Gamma(\Xm_{d}^{0}, \am^{\otimes h})^{G}.
\end{CD}
$$
\end{enumerate}
\end{proposition}

\begin{proof}[Proof]
Let $u \in G$. Since $\dim \am_{u \cdot F} = 1$, there is $c(u,F) \in k^{\times}$ such that
$$
(\varphi_{u}^{*} \eta)(F,\eta)=  c(u , F) \cdot \eta(F,\eta)= c(u , F).
$$
and $c$ is a ``crossed character", satisfying
$$
c(u u',F) =  c(u, F) \, c(u', u \cdot F ).
$$
Now the regular function $F \mapsto c(u, F)$ does not vanishes on
$\Xm_{d}^{0}$. By Lemma \ref{InvSections} below and the irreducibility of the
discriminant (Prop. \ref{discriminant}), we have
$$
c(u, F) = \chi(u) (\Disc F)^{n(u)}
$$
with $\chi(u)\in k^{\times}$ and $n(u) \in \ZZ$. The group $G$ being
connected, the function $n(u)=n$ is constant. Since
$c(\bfI_{3}, F) = 1$, we have $(\Disc F)^{n} =
\chi(\bfI_{3})^{- 1}$, and this implies $n = 0$. Hence, $c(u, F)$ is
independent of $F$ and $\chi$ is a character of $G$. Since the group of commutators of $G$ is $\SL_{3}(k)$, we have
$$
\chi(u) = \det(u)^{ w_0}
$$
for some $w_0 \in \ZZ$. It therefore suffices to calculate $\chi(u)$ when $u
= \lambda \bfI_{3}$, with $\lambda \in k^{\times}$. In this case $u \cdot F
= \lambda^{-d} F$. Moreover, the section $\eta_{f}$ is homogeneous of degree
$- 1$: if $\lambda \in k^{\times}$ and $F \in \Xm_{d}^{0}$, then
$$
\label{EtaHomo}
\eta_{f}(\lambda^{-d} F) / \eta_{f}(F)=  \lambda^{ d},
$$
hence,
$$
(\varphi_{u}^{*} \eta)(F,\eta)= \lambda^{ dg} = \det(u)^{w_0}.
$$
This implies
$$\lambda^{ 3 w_0} = \det(u)^{ w_0} = \lambda^{ dg},$$
and we have proven \eqref{SER2}.

Let $\Phi \in \Inv_{gh}(\Xm_{d}^{0})$ and $s = \tau(\Phi) = \Phi \cdot
\eta^{\otimes h}$, let also $w = d g h /3.$ Then
\begin{eqnarray*}
(\phi^*_u s)(F ,\eta) & = & \Phi(u \cdot F ) \cdot (\varphi^*_u \eta)(F,\eta)^{h}
\\ & = &
\det(u)^{-d (gh)/3} \Phi(F)  \cdot \det(u)^{ w_0 h}
\\ & = &
\Phi(F) = s(F,\eta),
\end{eqnarray*}
hence, $\tau(\Phi) \in \Gamma(\Xm_{d}^{0}, \lm^{\otimes
h})^{G}$. Conversely, the inverse of $\tau$ is the map $s \mapsto s/ \eta^{\otimes h}$, and this proves \eqref{SER3}.
\end{proof}

We made use of the following elementary lemma:

\begin{lemma}
\label{InvSections}
Let $f\in k[T_{1}, \dots, T_{n}]$ be irreducible and let $g \in k(T_{1},
\dots, T_{n})$ be a rational function which has neither zeroes nor poles
outside the set of zeroes of $f.$ Then there is an $m\in \ZZ$ and
$c\in k^{\times}$ such that $g=cf^m.$
\end{lemma}

\begin{proof}[Proof]
This is an immediate consequence of Hilbert's Nullstellensatz,
together with the fact that the ring $k[T_{1},\dots, T_{n}]$ is factorial.
\end{proof}

\subsection{Modular forms as invariants}
\label{PlaneCurves}

Let $d > 2$ be an integer and $g = \binom{d}{2}$. Since the fibres of
$\Ym_{d}^{0} \longrightarrow \Xm_{d}^{0}$ are  nonsingular non hyperelliptic
plane curves of genus $g$, by the universal property of $\Mm_g$ we get a
morphism
$$
\begin{CD}
p : \Xm_{g}^{0} & @>>> & \Mm_{g}^{0},
\end{CD}
$$
where $\Mm_{g}^{0}$ is the moduli stack of nonhyperelliptic curves of
genus $g$ and $p^{*}\lm = \am$ by construction. This induces a morphism
$$
\begin{CD}
p^{*}:
\Gamma(\Mm_{g}^{0}, \lm^{\otimes h})
& @>>> &
\Gamma(\Xm_{d}^{0}, \am^{\otimes h}).
\end{CD}
$$
Let $s \in \Gamma(\Mm_{g}^{0}, \lm^{\otimes h})$. For every $F \in
\Xm_{d}^{0}$, an element $u \in G$ induces an isomorphism
$$
\begin{CD}
\varphi_{u} : C_{ F } & @>>> & C_{u \cdot F}.
\end{CD}
$$
By the universal property of $\Mm_{g}^{0}$, the diagram
$$
\begin{CD}
\lm_{|p(X_d^0)} & @>{\Id}>> & \lm_{|p(X_d^0)} \\
@V{p^*}VV &    & @V{p^*}VV \\
\am & @>{\varphi^*_u }>> & \am
\end{CD}
$$
is commutative. Hence
$$\varphi_{u}^{*} \circ p^{*}(s) = p^{*}(s),$$
and this means that $p^*s \in \Gamma(\Xm_{d}^{0}, \am^{\otimes h})^G$.
Combining this result with Prop.\ref{Serre}\eqref{SER3}, we obtain:

\begin{proposition}
\label{DMToInv}
For any integer $h \geq 0$, the linear map $\sigma = {\tau^{- 1}} \circ p^{*}$
is a homomorphism:
$$
\begin{CD}
\Gamma(\Mm_{g}^{0}, \lm^{\otimes h})
& @>>> &
\Inv_{gh}(\Xm_{d}^{0})
\end{CD}
$$
such that
$$\sigma(f)(F) = f(C_{F},(p^*)^{-1} \eta)$$
for any $F \in \Xm_{d}^{0}$ and any section $f \in \Gamma(\Mm_{g}^{0},
\lm^{\otimes h})$.\qed
\end{proposition}

We now make a link between invariants and Siegel modular forms. Let $F \in X_d^0$ and
let $\eta_{1}, \dots, \eta_{g}$ be the basis of regular differentials on
$C_F$ defined in Sec.\ref{Quartics}. Let $\gamma_{1}, \dots \gamma_{2g}$ be a
symplectic basis of $H_{1}(C,\ZZ)$ (for the intersection pairing). The matrix
$$
\Omega = [\Omega_{1} \ \Omega_{2}] =
\begin{pmatrix}
\int_{\gamma_{1}}\eta_{1} & \cdots & \int_{\gamma_{2g}}\eta_{1}\\
\vdots &  & \vdots \\
\int_{\gamma_{1}}\eta_{g} & \cdots & \int_{\gamma_{2g}}\eta_{g}\\
\end{pmatrix}
$$
belongs to the set $\EuScript{R}_{g} \subset \bfM_{g,2g}(\CC)$ of
Riemann matrices, and $\tau = \Omega_{2}^{- 1} \Omega_{1} \in \HH_{g}$.
\begin{corollary}
\label{KleinGene1}
Let $f \in \bfS_{g,h}(\CC)$ be a geometric Siegel modular form,
$\widetilde{f} \in \bfR_{g,h}(\CC)$ the corresponding
analytic modular form, and $\Phi = \sigma( \theta^{*} f)$ the corresponding
invariant. In the above notation,
$$
\Phi(F) = (2 i \pi)^{gh}
\frac{\widetilde{f}(\tau)}{\det \Omega_{2}^{h}} \, .
$$
\end{corollary}

\begin{proof}
Let $\lambda = (p^{*})^{-1}(\eta)$ and $\omega = (\theta^{*})^{-1}(\lambda)$.
From Prop.\ref{Tau} and \ref{DMToInv}, we deduce
$$\Phi(F) = (\theta^{*}f)(C_{F},\lambda) = f(\Jac C,\omega).$$
On the other hand, Prop.\ref{Riemann} implies
$$
f(\Jac C ,\omega) = (2i \pi)^{gh}
\frac{\widetilde{f}(\tau)}{\det \Omega_{2}^{h}} \, ,
$$
from which the result follows.
\end{proof}

\section{The case of genus $3$}
\label{part3}

\subsection{Klein's formula}

We recall the definition of theta functions with (entire)
characteristics
$$
[\epsm] = \Ch{\varepsilon_{1}}{\varepsilon_{2}} \in \ZZ^g \oplus \ZZ^g,
$$
following \cite{lange}. The \emph{(classical) theta function} is given, for
$\tau \in \HH_{g}$ and $z \in \CC^g$, by
$$
\theta  \Ch{\varepsilon_{1}}{\varepsilon_{2}}(z, \tau) =
\sum_{n \in \ZZ^g}
q^{(n + \varepsilon_{1}/2) \tau (n + \varepsilon_{1}/2) + 2
(n + \varepsilon_{1}/2)(z+\varepsilon_{2}/2)}.
$$
The \emph{Thetanullwerte} are the values at $z = 0$ of these functions, and
we write
$$
\theta[\epsm](\tau) =
\theta \Ch{\varepsilon_1}{\varepsilon_2}(\tau) =
\theta \Ch{\varepsilon_1}{\varepsilon_2}(0,\tau).
$$
Recall that a characteristic is \emph{even} if $\varepsilon_1.\varepsilon_2
\equiv 0 \pmod{2}$ and \emph{odd} otherwise. Let $S_g$ (resp. $U_g$) be the
set of even characteristics  with coefficients in
$\{0,1\}$.
For $g \geq 2$, we put $h = \card{S_g} / 2 = 2^{g - 2}(2^g + 1)$ and
$$\widetilde{\chi}_h(\tau)=(2i\pi)^{gh} \prod_{\epsm \in S_g} \theta [\epsm](\tau).$$

In his beautiful paper \cite{igusa2}, Igusa proves the following result
[\loccit, Lem. 10 and 11]. Denote by $\widetilde{\Sigma}_{140}$ the modular form defined
by the thirty-fifth elementary symmetric function of the eighth power of the
even Thetanullwerte. Recall that a principally polarized abelian variety $(A,a)$ is decomposable if it is a product of principally polarized abelian varieties of lower dimension, and indecomposable otherwise.

\begin{theorem}
\label{igusath}
If $g \geq 3$, then $\widetilde{\chi}_h(\tau) \in \bfR_{g,h}(\CC)$.
Moreover, If $g = 3$ and $\tau \in \HH_3$, then:
\begin{enumerate}
\item
$A_{\tau}$ is decomposable if $\widetilde{\chi}_{18}(\tau) =
\widetilde{\Sigma}_{140}(\tau) = 0$.
\item
$A_{\tau}$ is a hyperelliptic Jacobian if $\widetilde{\chi}_{18}(\tau) = 0$
and $\widetilde{\Sigma}_{140}(\tau) \neq 0$.
\item \label{nonzero}
$A_{\tau}$ is a non hyperelliptic Jacobian if $\widetilde{\chi}_{18}(\tau)
\neq 0$.
\qed
\end{enumerate}
\end{theorem}

Using Prop. \ref{HoloGeom}, we define the geometric modular form of weight $h$
$$
\chi_{h}(A_{\tau}) = (2 i \pi)^{gh} \, \widetilde{\chi}_{h}(\tau)
(dz_{1} \wedge \cdots \wedge dz_{g})^{\otimes h}.
$$
Then Ichikawa \cite{Ichi3}, \cite{Ichi4} proved that $\chi_{h} \in
\bfS_{g,h}(\QQ)$. For $g = 3$, one finds
$$
\chi_{18}(A_{\tau}) = - (2 \pi)^{54} \, \widetilde{\chi}_{18}(\tau)
(dz_{1} \wedge dz_{2} \wedge dz_{3})^{\otimes 18}.
$$
Now we are ready to give a proof of the following result
\cite[Eq.~118, p.~462]{klein}:

\begin{theorem}[Klein's formula]
\label{kleinformula}
Let $F$ be a plane quartic defined over $\CC$ such that $C_{F}$ is
nonsingular. Let  $\eta_{1}, \eta_{2}, \eta_{3}$ be the classical basis of
$\Omega^{1}[C_{F}]$ and $\gamma_{1}, \dots \gamma_{6}$ be a symplectic basis
of $H_{1}(C_F,\ZZ)$ for the intersection pairing. Let
$$
\Omega = [\Omega_{1} \ \Omega_{2}] =
\begin{pmatrix}
\int_{\gamma_{1}}\eta_{1} & \cdots & \int_{\gamma_{6}}\eta_{1}\\
\vdots &  & \vdots \\
\int_{\gamma_{1}}\eta_{3} & \cdots & \int_{\gamma_{6}}\eta_{3}\\
\end{pmatrix}
$$
be a period matrix of $\Jac(C)$ and $\tau = \Omega_{2}^{- 1} \Omega_{1} \in
\HH_3$. Then
$$
\Disc(F)^2 = \frac{1}{2^{28}}(2 \pi)^{54}
\frac{\widetilde{\chi}_{18}(\tau)}{\det(\Omega_2)^{18}}.
$$
\end{theorem}

\begin{proof}
The Cor.\ref{KleinGene1} shows that $I = \sigma \circ \theta^{*}(\chi_{18})$ is an invariant of weight $54$, and for any $F \in \Xm_4^0$,
$$
I(F) = - (2 \pi)^{54}
\frac{\widetilde{\chi}_{18}(\tau)}{\det \Omega_{2}^{18}}.
$$
Moreover Th. \ref{igusath}\eqref{nonzero} shows that $I(F) \ne 0$
for all $F \in \Xm_4^0$. Applying Lem. \ref{InvSections} for the
discriminant, we find by comparison of the weights that
$I = c \Disc^2$ with $c \in \CC$ a constant. But if $F_{m}$ is the Ciani
quartic associated to a matrix $m \in \Sym_{3}(k)$ as in Example
\ref{CianiQuartics}, and if $\Disc F_{m} \neq 0$, then it is proven in
\cite[Cor. 4.2]{SAGA} that Klein's formula is true for $F_{m}$ and $c = -
2^{28}$.
\end{proof}

\begin{remark}
\label{mu9}
The morphism $\theta^{*}$ defines an injective morphism of graded $k$-algebras
$$
\begin{CD}
\bfS_{3}(k) = \oplus_{h \geq 0} \bfS_{3,h}(k) & @>>> &
\bfT_{3}(k) = \oplus_{h \geq 0} \bfT_{3,h}(k).
\end{CD}
$$
In \cite{Ichi2}, Ichikawa proves that if $k$ is a field of characteristic
$0$, then $\bfT_{3}(k)$ is generated by the image of $\bfS_{3}(k)$ and
a primitive Teichm\"{u}ller form $\mu_{3,9}\in \bfT_{3, 9}(\ZZ)$ of weight $9$, which is not of Siegel
modular type. He also proves in \cite{Ichi4} that
\begin{equation}
\label{chimu}
\theta^{*}(\chi_{18}) = -2^{28} \cdot (\mu_{3,9})^2.
\end{equation}

Th. \ref{kleinformula} implies that $\mu_{3,9}$ is actually equal to the discriminant up to a sign. This might probably be deduced from the definition of $\mu_{3,9}$, although we did not sort it out (see also \cite[Sec. 2.4]{jong}).
\end{remark}

\begin{remark}
Besides  \cite{lockhart} and \cite{guardia} where an analogue of Klein's
formula is derived in the hyperelliptic case, there exists a beautiful
algebraic Klein's formula, linking the discriminant with
irrational invariants \cite[Th.11.1]{gizatullin}.
\end{remark}

\subsection{Jacobians among abelian threefolds}

Let $k\subset\CC$ be a field and let $g = 3$. We prove the following theorem
which allows to determine whether a given abelian threefold defined over $k$
is $k$-isomorphic to a Jacobian of a curve defined over the same field. This
settles the question of Serre recalled in the introduction.

\begin{theorem}
\label{principal}
Let $(A,a)$ be a principally polarized abelian threefold defined over $k
\subset \CC$.  Let $\omega_{1}, \omega_{2}, \omega_{3}$ be a basis of
$\Omega^{1}_{k}[A]$ and $\gamma_{1}, \dots \gamma_{6}$ a symplectic basis
of $H^1(A,\ZZ)$, in such a way that
$$
\Omega = [\Omega_{1} \ \Omega_{2}] =
\begin{pmatrix}
\int_{\gamma_{1}} \omega_{1} & \cdots & \int_{\gamma_{6}}\omega_{1}\\
\vdots &  & \vdots \\
\int_{\gamma_{1}}\omega_{3} & \cdots & \int_{\gamma_{6}}\omega_{3}\\
\end{pmatrix}
$$
is a period matrix of $(A, a)$. Put $\tau = \Omega_{2}^{- 1} \Omega_{1} \in
\HH_3$.
\begin{enumerate}
\item
If $\widetilde{\Sigma}_{140}(\tau)=0$ then   $(A,\lambda)$ is
decomposable. In particular it is not a Jacobian.
\item
If $\widetilde{\Sigma}_{140}(\tau) \ne 0$ and
$\widetilde{\chi}_{18}(\tau)=0$ then there exists a hyperelliptic curve $X/k$
such that  $(\Jac X,j) \simeq (A,a)$.
\item
If $\widetilde{\chi}_{18}(\tau) \ne 0$ then $(A,a)$ is isomorphic to a
Jacobian if and only if
$$-\chi_{18}(A,\omega_1 \wedge \omega_2 \wedge \omega_3)=(2 \pi)^{54} \frac{\widetilde{\chi}_{18}(\tau)}{\det(\Omega_2)^{18}}$$
is a square in $k$.
\end{enumerate}
\end{theorem}

\begin{proof}
The first and second points follow from Th.\ref{igusath} and Th.\ref{twist}.
Suppose now that $(A, a)$ is isomorphic over $k$ to the Jacobian of a non
hyperelliptic genus $3$ curve $C/k$. Let $F \in X_4^0$ be a plane model of
$C$. Using Th.\ref{kleinformula} we get that
$$
-\chi_{18}(A,\omega_1 \wedge \omega_2 \wedge \omega_3)=(2 \pi)^{54}
\frac{\widetilde{\chi}_{18}(\tau)}{\det(\Omega_2)^{18}}=2^{28} \Disc(F)^2
$$
so it is a square in $k$. On the contrary, Cor.\ref{ftwist} shows that if
$(A',a')$ is a quadratic twist of a Jacobian $(A,a)$ then the expression
$$-f(A',\omega_1' \wedge \omega_2' \wedge \omega_3')=(2 \pi)^{54} \frac{\widetilde{\chi}_{18}(\tau')}{\det(\Omega_2')^{18}}$$
is not a square.
\end{proof}

\begin{remark} \label{mu9serre}
Note that one does not really need  Klein's formula. Alternatively, we could use
\eqref{chimu} which also proves  that $\theta^*(-\chi_{18})$ is a square.
\end{remark}

\begin{corollary}
In the notation of Th.\ref{principal}, the quadratic character $\varepsilon$
of $\Gal(k_{\textrm{sep}}/k)$ introduced in Theorem \ref{twist}
is given by $\varepsilon(\sigma) = d/ d^{\sigma}$, where
$$
d = \sqrt{(2 \pi)^{54}
\frac{\widetilde{\chi}_{18}(\tau)}{\det(\Omega_{2})^{18}}} \, ,
$$
with an arbitrary choice of the square root.
\end{corollary}

\subsection{Beyond genus $3$} \label{beyond}

It is natural to try to extend our results to the case $g > 3$. The first question to ask is

\begin{itemize}
\item Does there exist an analogue of Klein's formula for $g > 3$ ?
\end{itemize}

Here we can give a partial answer. Using Sec.\ref{sec_Teichmuller}, we can
consider the Teichm\"{u}ller modular form $\theta^{*}(\chi_h)$ with $h
= 2^{g - 2} (2^g + 1)$. In \cite[Prop.4.5]{Ichi4} (see also \cite{Tsu}), it
is proven that for $g > 3$ the element

$$\theta^*(\chi_h)/2^{2^{g-1}(2^g-1)}$$

has as a square root a primitive element $\mu_{g,h/2} \in \bfT_{g,h/2}(\ZZ)$.
If $g=4$, in the footnote, p. 462 in \cite{klein} we find the following
amazing formula
\begin{equation}
\label{g4}
\frac{\widetilde{\chi}_{68}(\tau)}{\det(\Omega_2)^{68}}=c \cdot \Delta(X)^2
\cdot T(X)^8.
\end{equation}
Here $\tau=\Omega_{2}^{-1}\Omega_{1}$, with $\Omega=[\Omega_{1} \ \Omega_{2}]$ beeing a period matrix of a genus $4$ non hyperelliptic curve $X$ given in $\PP^3$ as an intersection of a quadric $Q$ and a cubic surface
$E$.  The elements $\Delta(X)$ and
$T(X)$ are defined in the classical invariant theory as, respectively, the
discriminant of $Q$ and the tact invariant of $Q$ and $E$ (see
\cite[p.122]{Salmon}). No such formula seems to be known in the non
hyperelliptic case for $g > 4$.

Let us now look at what happens when we try to apply Serre's apporoach for $g > 3.$ To begin with, when $g$ is even, we cannot use Cor.\ref{ftwist} to distinguish between quadratic twists. In particular, using the previous result, we see that
$\chi_h(A,\omega_k)$ is a square when $A$ is a principally polarized abelian
variety defined over $k$ which is geometrically a Jacobian. A natural
question is:

\begin{itemize}
\item
What is the relation between this condition and the locus of geometric
Jacobians over $k$?
\end{itemize}

Let us assume now that $g$ is odd. As we pointed out in Rem.\ref{mu9serre},
the existence of the square root is almost sufficient to answer Serre's
questions when $g=3$.  This is not the case when $g>3$.  The proof of the corollary
\ref{cortwist} shows that
$$\overline{\chi}_h(A')=(d^2)^{h/2} \overline{\chi}_h(A)$$
for a Jacobian $A$ and a quadratic twist $A'$.
 What enables us to distinguish between $A$ and $A'$ when $g=3$ is the
following: if $A$ is the Jacobian of a curve then  $\chi_h(A)$ is a square
whereas $d^2$ is not and $h/2=9$ is odd. However as soon as $g>3$, $2\mid
2^{g-3}$, the power $g-3$ being the maximal power of 2 dividing $h/2$, so it
is not enough for $\widetilde{\chi}(A)$ to be a square in $k$ to make a
distinction between $A$ and $A'$. It must rather be an element of
$k^{2^{g-2}}$.

It can be easily seen from  the proof of  \cite[Th.1]{Tsu} that $\theta^*(\chi_h)$ does not admit a fourth
root. According to \cite{Bilu} or \cite{Xarles} this implies
 $\overline{\chi}_h(A) \notin k^{2^{g-2}}$ for infinitely many
 Jacobians $A$ defined over number fields $k$.
 So we can no longer use the modular form $\chi_h$ to
 easily characterize
  Jacobians over $k$. So the question is:

\begin{itemize}
\item
Is it possible to find a modular form to replace $\chi_h$ in Serre's strategy when $g>3$ ?
\end{itemize}

\UnTrait
\end{document}